\documentclass[12pt,a4paper]{article}
\usepackage{latexsym,amsmath,amssymb}
\usepackage{epsfig,array}

\def\reff#1{(\ref{#1})}
\def\T{^\top}

\newtheorem{thm}{Theorem}[section]
\newtheorem{prop}{Proposition}[section]
\newenvironment{proof}{\noindent {\bf Proof. }}{\hfill $\Box$ \newline\par}
\def\reff#1{(\ref{#1})}
\def\qf#1#2{#1\T  #2 #1}

\def\T{^\top}
\newcommand{\R}{\mathbb R}

\def\PP{{\mathcal P}}
\def\NN{{\mathcal N}}

\def\CC{{\mathcal C}}
\def\Ko{{\mathcal K}_0}
\def\Kos{{\mathcal K}_0^*}

\def\Diag{\mbox{Diag }}
\def\diag{\mbox{diag }}

\def\beq#1{\begin{equation}\label{#1}}
\def\eeq{\end{equation}}
\def\bep{\begin{proof}}
\def\ep{\end{proof}}

\def\lk{\left\{}
\def\rk{\right \} }
\def\norm#1{\|#1\|}

\begin{document}
\begin{titlepage}
{\Large \raggedright Improving SDP bounds for minimizing quadratic functions over the $\ell^1$-ball
\par
}
\begin{flushleft}
  \normalsize
Immanuel M. Bomze\\
Florian Frommlet\\
Martin Rubey\\ [2ex] all: ISDS, University of Vienna \\[1.2ex]
Mail address: Universit\"atsstra{\ss}e 5/9, A-1010 Wien, Austria
\end{flushleft}
\rule{\textwidth}{0.75pt}

\begin{abstract}
  \normalsize
\noindent  In this note, we establish superiority of the so-called copositive bound over
  a bound suggested by Nesterov for the problem to minimize a
  quadratic form over the $\ell ^1$-ball. We illustrate the improvement by
  simulation results using Jos Sturm's SeDuMi. The copositive bound has the additional advantage that it
  can be easily extended to the inhomogeneous case of quadratic objectives
  including a linear term. We also indicate some improvements of the eigenvalue bound for the quadratic optimization over the $\ell ^p$-ball with
  $1<p<2$, at least for $p$ close to one.
\end{abstract}

\begin{flushleft}
  \normalsize

\noindent

This version: \today
\end{flushleft}
\rule{\textwidth}{0.75pt}

\hbox{}\hfill{\em This work is dedicated to the memory of Jos Sturm.}
\end{titlepage}
\newpage
\setcounter{page}{1}
\section{\label{intro}Introduction}
As usual we denote for $p>0$ by ${\norm y}_p = \left [\sum_{i=1}^n |y_i|^p \right ]^{1/p}$
the $p$-norm of a vector $y$ in $n$-dimensional Euclidean space $\R^n$. The $\ell^p $-ball is then $B_p=\lk y\in \R ^n :  {\norm y}_p \le 1\rk$.
In the Handbook of
Semidefinite Programming \cite{NestWolkYe:00}, Nesterov deals with the problem of minimizing a quadratic form
over the $\ell^1$-ball $B_1$, among many others (see p.~387 in \cite{NestWolkYe:00}), and specifies a
bound obtained by SDP relaxation for this problem. However, the feasible set $B_1$ is a polytope with not too many, and known,
vertices, namely $\pm e_i$, $1\le i\le n$, where $e_i$ are vectors in $\R^n$ with one as $i$-th coordinate
and zeroes elsewhere. This elementary observation allows for transformation of the problem into a moderately
sized Standard Quadratic Optimization Problem (StQP), which consists of minimizing a quadratic form
over the standard simplex
$$\Delta^{2n} = \lk x\in \R^{2n}_+ : \bar e\T x = 1\rk\, ,$$
where $\bar e = [e\T , e\T ]\T\in \R^{2n}$ and $e=\sum_i e_i \in \R ^n$ are all-ones vectors in $\R^{2n}$ and
$\R^n$, respectively ($\R ^m_+$ denotes the set of all vectors in $\R ^m$ with no negative coordinates).
This way, we may apply any valid bound for StQPs also to 
the quadratic problem over $B_1$. Here, we focus on the so-called copositive (relaxation) bound introduced in \cite{BdeK}
and recently investigated in more detail in \cite{AnBu}.

The paper is organized as follows: First, we present the reformulation of the problem as a StQP, then
investigate the copositive bound for this special problem class and establish superiority compared to the
bound proposed by Nesterov.  We also address the special case of sign-constrained data matrices where it can
be shown that the size of the StQP can be kept to the original size, so that -- unlike the general case --
doubling the dimension can be avoided. Empirical quality assessment is provided by a small simulation study.
Finally, we use the established copositive bound to obtain a lower bound when minimizing a quadratic form
over the $\ell^p$-ball $B_p$ for $1<p<2$.

\section{\label{transf} Reformulation as an StQP} 

Given a symmetric, possibly indefinite $n\times n$ matrix $C$ and a vector $c\in \R^n$, consider the
quadratic minimization problem over the $\ell^1$-ball
\beq{el1def} \alpha ^* = \min \lk \qf yC + 2c\T y : y\in B_1\rk \, . \eeq
Since the feasible set $B_1$ of the above problem is the convex hull of the vectors $\pm
e_i$ introduced above, we can write
\beq{transstqp} \alpha ^* = \min \lk  x\T (V\T CV+ \bar e c\T V + V\T c
\bar e)x : x\in \Delta ^{2n}\rk\, ,
\eeq
where $V= [I_n , -I_n ]$ contains all vertices of the $\ell^1$-ball
as column vectors.  Hence, introducing the symmetric $(2n)\times (2n)$ matrix
$$Q_{C,c} = \left[ \begin{array}{rr}
      C  & -C \\
      -C &  C
   \end{array}\right] +
   \left[ \begin{array}{rr}
     e c\T + ce\T   & ce\T-e c\T \\
       e c\T -ce\T   & -e c\T - ce\T
   \end{array}\right]
$$
we arrive at the StQP
$$\alpha^* = \min \lk \qf x{Q_{C,c}} : x\in \Delta ^{2n}\rk\, ,$$
and any StQP bound for $Q_{C,c}$ is a valid bound for $\alpha ^*$. Here, we will focus on the so-called
copositive relaxation bound. For the readers' convenience, we provide some background in the following
section.

\section{\label{coprelbd} Copositive relaxation bounds for StQPs}

As is well known~\cite{BDdeKQRT}, we can reformulate every StQP in $m$-dimensional space
of the form
\beq{stqp} \alpha_Q^* = \min \lk
\qf x Q : x\in \Delta^m \rk \eeq
into a copositive program:
\beq{cop} \max \left\{ \lambda : Q - \lambda E
\in \CC \right \} = \alpha_Q^*\, , \eeq
where $E$ is the all-ones $m\times m$ matrix and
\beq{copcon} \CC =
\lk A \in \mathcal S:  \qf xA \ge 0 \mbox{ for all } x\in \R^m_+ \rk \eeq
denotes the convex cone of all {\em
copositive} matrices. Here $\mathcal S$ are the symmetric $m\times m$ matrices.

Let now $\PP \subset \mathcal S$ denote the cone of all positive semidefinite matrices and $\NN \subset
\mathcal S$  the cone of all nonnegative symmetric matrices. Then, a (zero-order) approximation \cite{BdeK}
of the copositive cone is given by $\Ko = \PP + \NN \subseteq \CC$ with $\CC = \Ko$ only if $n\le 4$ \cite{Diananda:67}.
Replacing $\CC$ with $\Ko$ yields the {\em copositive relaxation bound:}
\beq{alf1}
 \alpha _Q ^{\rm cop} = \max \left\{ \lambda : Q - \lambda E \in \Ko = \PP + \NN \right \} \le \alpha^* _Q\, .
\eeq

Passing to the dual problems of~\reff{cop} and \reff{alf1}, we obtain alternative formulations for
$\alpha^*_Q$ and $\alpha _Q^{\rm cop}$, respectively.  To this end, we need some more notation. Let $A,B$ be
two symmetric $m \times m$ matrices and recall that the \emph{trace} of a matrix is the sum of its diagonal
elements. Then define the {\em inner product}
\[
A \bullet B = \mbox{trace}(AB) = \sum_{ij}a_{ij}b_{ij}\, .
\]
Now strong duality arguments yield equality of~\reff{cpp} below with~\reff{cop}, and~\reff{alf1}
with~\reff{alf1dual} below, respectively. For details see~\cite{BdeK}.
\beq{cpp} \min  \left\{ Q\bullet X :
E\bullet X = 1 \, , \, X\in \CC^* \right \} = \alpha^* _Q\, , \eeq
with $\CC^*$ the completely positive cone,
the dual cone of $\CC$ given by
\beq{comcon} \CC^* = \left \{\sum_{i=1} ^k y_i (y_i)\T  \in \mathcal{S}: y_i\in \R ^m _+ \, , \, i\in\lk 1, \ldots , k\rk\, , \hbox { some }
k  \right \}\, ,
\eeq
and
\beq{alf1dual} \min  \left\{ Q\bullet X : E\bullet X = 1 \, , \, X\in \PP \cap \NN \right \}= \alpha
_Q ^{\rm cop} \, , \eeq as $\Ko= \PP + \NN  $ has the dual cone $\Kos = \PP \cap \NN \supset \CC^*$.

A direct argument why the solution of~\reff{alf1dual} never can exceed $\alpha ^* _Q$ employs the fact that
for any $x\in \Delta^m$, the rank-one matrix $X= xx\T$ satisfies $X\in \Kos$ as well as $E\bullet X = (e\T
x)^2 = 1$, along with $Q\bullet X= \qf x Q$.

\section{\label{copel1}Copositive relaxation bounds\\
for the $\ell ^1$-constrained problem}

Now let us specialize to $Q=Q_{C,c}$, where $m=2n$. We focus on the dual formulation~\reff{alf1dual}, and
partition any $X\in\PP \cap \NN $ in a natural way:
$$ X= \left[ \begin{array}{rr}
      U  & Y \\
      Y\T &  V
    \end{array}\right]\, .
$$
Straightforward calculations then yield
$$Q_{C,c}\bullet X = C\bullet (U-Y-Y\T +V) + (ec\T+ce\T)\bullet (U-V) - (ec\T-ce\T)\bullet (Y-Y\T)\, .$$
In order to compare the copositive bound with the bound suggested by Nesterov we will restrict ourselves to
the homogeneous case $c=o$, and abbreviate $Q_C = Q_{C,o}$\,. Then the copositive bound is given by
\beq{el1dual} \alpha _{Q_C} ^{\rm cop} = \min \lk C\bullet (U-Y-Y\T +V) : X= \left[ \begin{array}{rr}
      U  & Y \\
      Y\T &  V
    \end{array}\right] \in \PP \cap \NN \, ,\;  E\bullet X = 1\right \}\, .
\eeq

On p.~387 of \cite{NestWolkYe:00}, Nesterov proposes an alternative SDP relaxation bound as follows (in
maximization rather than minimization form):
\beq{el1nest} \alpha _{C}
^{\rm Nest} = \min \lk C \bullet W : \Diag (u)\succeq W \succeq O \mbox{ for some }u\in \Delta ^n\rk \eeq
(actually he requires $e\T u \le 1$ rather than $e\T u = 1$ as included in the definition of $\Delta^n$, but
obviously, scaling $u$ such that $e\T u = 1$ holds does not change \reff{el1nest}-feasibility). Here and in
the sequel, we abbreviate by $A\succeq B$ the fact that $A-B$ is positive-semidefinite.

Our aim is now to show that always $\alpha _{C} ^{\rm Nest}\le \alpha _{Q_C} ^{\rm cop} \le \alpha ^*$ holds,
where the last inequality follows from validity of the general StQP bound. Note that in the (easy convex)
case of positive semidefinite $C$, all three values are equal to 0. In order to prove the first inequality in the
general case, we start with an auxiliary result.

\newtheorem{lem}{Lemma}

\begin{lem}\label{matree}
  Let
  $$
      X=
    \left[ \begin{array}{cc}
         U   &  Y \\
       Y\T &    V
     \end{array}\right]
  $$
  be a symmetric positive-semidefinite matrix where all sub-matrices are square
  and of the same size, and $ U$, $V$ and $Y$ have no negative entries. Define
  $$
      X_-=
  \left[ \begin{array}{cc}
        U   &   -Y \\
        -Y\T  &   V
    \end{array}\right]
  $$
  and let $D= \Diag (X\bar e)$ be the diagonal matrix containing the row-sums
  of $X$ on the diagonal. Then
  \begin{description}
  \item[(a)] $X_-$  is positive-semidefinite; and
  \item[(b)] $  D-  X_-$ is positive-semidefinite.
  \end{description}
\end{lem}
\begin{proof}
  First we show that $X\in \PP$  if and only if $X_- \in \PP$. Because of the special structure of $ X$ and $ X_-$,
  we have that for any two vectors $a$ and $b$
  \beq{swapdef}
   \left[ \begin{array}{rr}
       a\T  & b\T
    \end{array}\right]
      X_-
   \left[ \begin{array}{c}
       a \\ b
    \end{array}\right]
    = a\T    U a - 2 a\T    Y b + b\T    V b
    =
    \left[ \begin{array}{rr}
       a\T  & -b\T
     \end{array}\right]
      X
   \left[ \begin{array}{r}
       a \\ -b
    \end{array}\right]\, .
  \eeq
  Hence assertion~(a) is established.
  Furthermore, since  $D$ is a diagonal
  matrix, it has only zero off-diagonal blocks.
     Therefore, by the same argument as before, $ D-X\in \PP$
  if and only if $ D-X_-\in \PP$, so it remains to show that $D- X$ is positive-semidefinite.

  To this end we interpret $X$ as the adjacency matrix of an edge-weighted graph, and $D-X$ as its corresponding
  Laplacian matrix, which is well known to be positive-semidefinite, if there are no negative weights. For the convenience of the reader we give a
  short proof, for more details see e.g. \cite{Biggs}. The Laplacian $D-X$ can be rewritten as $BGB\T$, where $B$ is
  the oriented incidence matrix and $G$ is a diagonal matrix containing the weights of the graph. Indeed,
  indexing the columns of the $m\times\binom{m}{2}$
  matrix $B$ by ordered pairs $(k,l)$ with $1\leq k<l\leq m$, we have
  $$
    B_{i,(k,l)}=
    \begin{cases}
      +1&\text{ if $i=k$}\\
      -1&\text{ if $i=l$}\\
      \phantom{-}0&\text{ otherwise}
    \end{cases}
  $$
  and correspondingly
  $$
    G_{(k,l),(r,s)}=
    \begin{cases}
      X_{(k,l)}&\text{ if $(k,l)=(r,s)$}\\
      0        &\text{ otherwise.}
    \end{cases}
  $$
It is then straightforward to show that $BGB\T = D-X$. Now, since $X$ has no negative entries, we can take
the square root $G^{\frac{1}{2}}$  of the diagonal matrix $G$ and thus obtain $D-X =(B
G^{\frac{1}{2}})(BG^{\frac{1}{2}})\T \succeq O$.
\end{proof}

\begin{thm}
We have $\alpha _{C} ^{\rm Nest}\le \alpha _{Q_C} ^{\rm cop} $ for all symmetric $n\times n$ matrices $C$,
with strict inequality for some instances.
\end{thm}

\begin{proof}
  We show that for any \reff{el1dual}-feasible $X$, the matrix $W=U-Y- Y\T + V$
  together with the vector $u=(U+Y+Y\T + V)e$ satisfies the constraints of
  \reff{el1nest}. Then the assertion follows immediately by definition of the
  problems defining the two bounds.\\
  First, Lemma~\ref{matree}(a) ensures
  positive semidefiniteness of $W$, since $\qf yW = [y\T, y\T ]X_-[y\T ,y\T ]\T
  \ge 0$. Next, $u\in \R^n_+$ by the nonnegativity assumption on $U$, $Y$, and
  $V$, and $e\T u = \bar e\T X \bar e = E\bullet X = 1$, whence $u\in \Delta
  ^n$ results. Further, $D= \Diag (X\bar e)$ satisfies
$$\qf y{(\Diag u)}= 
\sum_i [(U+Y)e]_i y_i^2 +  \sum_i [(Y\T +V)e]_i y_i^2 = [y\T, y\T ]D\left[ \begin{array}{r}
       y \\ y
    \end{array}\right]\, ,$$
whence $\qf y{(\Diag u - W)}=[y\T, y\T ](D-X_-)[y\T ,y\T ]\T \ge 0$ for all $y\in \R^n$ results, due to
Lemma~\ref{matree}(b). Therefore we conclude $\Diag u \succeq W\succeq O$, and $(W,u)$ is
\reff{el1nest}-feasible. Finally, to establish strict domination of the copositive bound, it suffices to look at the $3\times 3$ instance
\beq{strictex}
C=\left[ \begin{array}{rrr} -1 &a &-1 \\
a &-1 &-1 \\-1 &-1 &-1
\end{array}\right] \, .
\eeq
For $a=1$ we get $\alpha _{C} ^{\rm Nest} = -4/3$ whereas $\alpha _{Q_C} ^{\rm cop}= -1 = \alpha^* $ is
even exact.
\end{proof}

Thus, by using a StQP reformulation of the original problem we obtain a better SDP bound. In case one does not want to solve the SDP relaxations to optimality and
still needs a valid bound, the dual formulations could be preferred, which give a valid bound for any feasible objective value.
For completeness, we specify them for both bounds:
$$\alpha _{Q_C} ^{\rm cop} = \max\lk \lambda \in \R : \lambda \le (Q_C)_{ij} - \bar S_{ij} \mbox{ for all }i,j\, ,\mbox{ and some }\bar S\succeq O\rk$$
which is the originally primal form~\reff{alf1},
and the dual formulation of the Nesterov bound which with small modifications appears already in~\cite{NestWolkYe:00}:
$$\alpha _{C} ^{\rm Nest} = \max \lk  \lambda \in \R : \diag (S) \le -\lambda e\, ,\, S\succeq - C \mbox{ for some }  S\succeq O\rk \, .$$
A direct comparison of these two formulations does not seem immediately evident.

The price we have to pay
for more efficiency of $\alpha _{Q_C} ^{\rm cop}$ is the double dimension. One may wonder what happens if the $(2n)\times (2n)$ matrix
$Q_C$ is replaced by $C$ itself, thus arriving at an SDP of the same size as that used in $\alpha _{C} ^{\rm
Nest}$.

\begin{thm}\label{smallcop}
For any symmetric $C$, we have
$$\alpha _{Q_C} ^{\rm cop} \le \alpha _{C} ^{\rm cop}\, ,$$
but in general, the latter also exceeds the true value $\alpha^*$.
\end{thm}

\begin{proof}
Decompose any \reff{el1dual}-feasible $X$ again as
$$X=
    \left[ \begin{array}{cc}
         U   &  Y \\
       Y\T &    V
     \end{array}\right]\, .
 $$
If, in particular, $Y=V=O$, we arrive at $C\bullet (U-2Y+V) = C\bullet U$, and the conditions $X\in \PP \cap
\NN$ and $E\bullet X = 1$ boil down to $U\in \PP \cap \NN$ (by abuse of notation, we ignore differences of
matrix size in the cones) as well as $ee\T \bullet U = 1$, so that
$$\alpha _{Q_C} ^{\rm cop} \le \min \lk C\bullet U : ee\T \bullet U = 1\, , \; U\in \PP \cap \NN\rk = \alpha _{C} ^{\rm cop}\, ,$$
which proves the assertion. For the matrix $C$ of (\ref{strictex}) with $a=2$ we have $-1.5 = \alpha^*
<\alpha _{C} ^{\rm cop}=-1$.
\end{proof}

\section{Sign-constrained data matrices}

Here we focus on the particular case where $C$ has no positive entries, i.e. where $-C\in \NN$. In this case
the minimum of the quadratic form over the $\ell ^1$-ball is attained at a point of the standard simplex:
\begin{prop}\label{negprop}
For any symmetric $C$ with no positive entries we have
$$\alpha^* = \min \lk \qf xC : x\in B_1\rk  = \min \lk \qf xC : x\in \Delta\rk \, .$$
\end{prop}
\begin{proof}
First note that except for the trivial case $C=O$ a nonnegative matrix cannot be positive semidefinite.
Therefore it is guaranteed that $\alpha^*$ is attained at a point $\tilde x$ on the boundary of $B_1$, whence we get
$\sum_{i=1}^n |\tilde x_i| = 1$. Then  $x^*:=\left [ |\tilde
x_1|, \ldots , |\tilde x_n|\right]\T \in \Delta ^n$ and from $-C\in \NN$ it follows $(x^*)\T Cx^* \le (\tilde x)\T C \tilde x$, and
thus the proposition.
\end{proof}
This relation has an exact counterpart for the SDP-relaxation: we can indeed use the cheaper copositive bound $\alpha
_{C} ^{\rm cop}$ considered at the end of the previous section, because the latter coincides with $\alpha
_{Q_C} ^{\rm cop} $.

\begin{thm}\label{negc}
For any symmetric $C$ with no positive entries, we have
$$\alpha _{Q_C} ^{\rm cop} = \alpha _{C} ^{\rm cop}\le \alpha^*\, ,$$
and again, in general, $\alpha _{C} ^{\rm cop} > \alpha _{C} ^{\rm Nest}$ for some instances $C\in -\NN$.
\end{thm}
\begin{proof}
In view of Theorem~\ref{smallcop}, we only have to show $\alpha _{Q_C} ^{\rm cop} \ge \alpha _{C} ^{\rm
cop}$. Now, if
$$X=
    \left[ \begin{array}{cc}
         U   &  Y \\
       Y\T &    V
     \end{array}\right]\, \in \PP \cap \NN\, ,
 $$
and also $-C\in \NN$, we get $C\bullet (U-Y-Y\T + V) \ge C\bullet (U+Y+Y\T + V) $. Now $Z:=U+Y + Y\T + V$ is
a symmetric $n\times n$ matrix which has no negative entries and is also positive-semidefinite, as follows
from
$$\qf yZ = [y\T , y\T ] X [y\T , y\T ]\T \ge 0\quad\hbox{ for all }y\in \R ^n\, .$$
Further, $ee\T\bullet Z = E\bullet X$, which shows \beq{negceq}
\begin{array}{rl}
\alpha _{Q_C} ^{\rm cop} &= \min \lk C\bullet (U-2Y +V) : X= \left[ \begin{array}{rr}
      U  & Y \\
      Y &  V
    \end{array}\right] \in \PP \cap \NN \, ,\;  E\bullet X = 1\right \} \\
&\ge \min \lk C\bullet Z : Z \in \PP \cap \NN \, ,\;  ee\T\bullet Z = 1\right \}= \alpha _{C} ^{\rm cop} \, .
\end{array} \eeq Obviously the matrix $C$ of (\ref{strictex}) with $a=0$ has no positive entries and we
have $\alpha _{C} ^{\rm Nest} \approx -1.1429$ whereas $\alpha _{Q_C} ^{\rm cop}= \alpha^* = -1$ is
exact.\end{proof}

\section{\label{emp} Empirical findings}

To compare the quality of the two SDP relaxation bounds we generated 1000 random symmetric $n \times n$
matrices for $n=10$ and $n=20$ respectively. The SDP problems were solved using Jos Sturm's SeDuMi \cite{Sturm:99}. To
obtain upper bounds we applied two optimization procedures for the StQP \reff{transstqp}, a fixed-step
exponential replicator dynamics with $\theta = 0.05$, and Wolfe's reduced gradient method. For details and
comparison with other local optimization procedures for StQPs see \cite{portfo}. For both iterative methods
we used 10 random starting vectors respectively and took the overall minimum of both procedures as a
reference solution $\alpha  ^{\rm ref}$. We thus have $ \alpha _{C} ^{\rm Nest}\le \alpha _{Q_C} ^{\rm cop}
\le \alpha  ^{\rm ref}$ and in case of equality of the last two bounds we can conclude $\alpha ^{\rm ref} =
\alpha _{C} ^{\rm cop} = \alpha^*$.

\begin{figure} 
\begin{center}
\caption{Distribution of  $ \alpha _{Q_C} ^{\rm cop}/ \alpha _{C} ^{\rm Nest}$ for 1000 randomly generated
$n\times n$ matrices. } $\begin{array}{c@{\hspace{1in}}c} \multicolumn{1}{l}{\mbox{\bf (a)} \quad n=10} &
    \multicolumn{1}{l}{\mbox{\bf (b)} \quad n=20} \\ [-0.03cm]
\epsfxsize=2.6in \epsffile{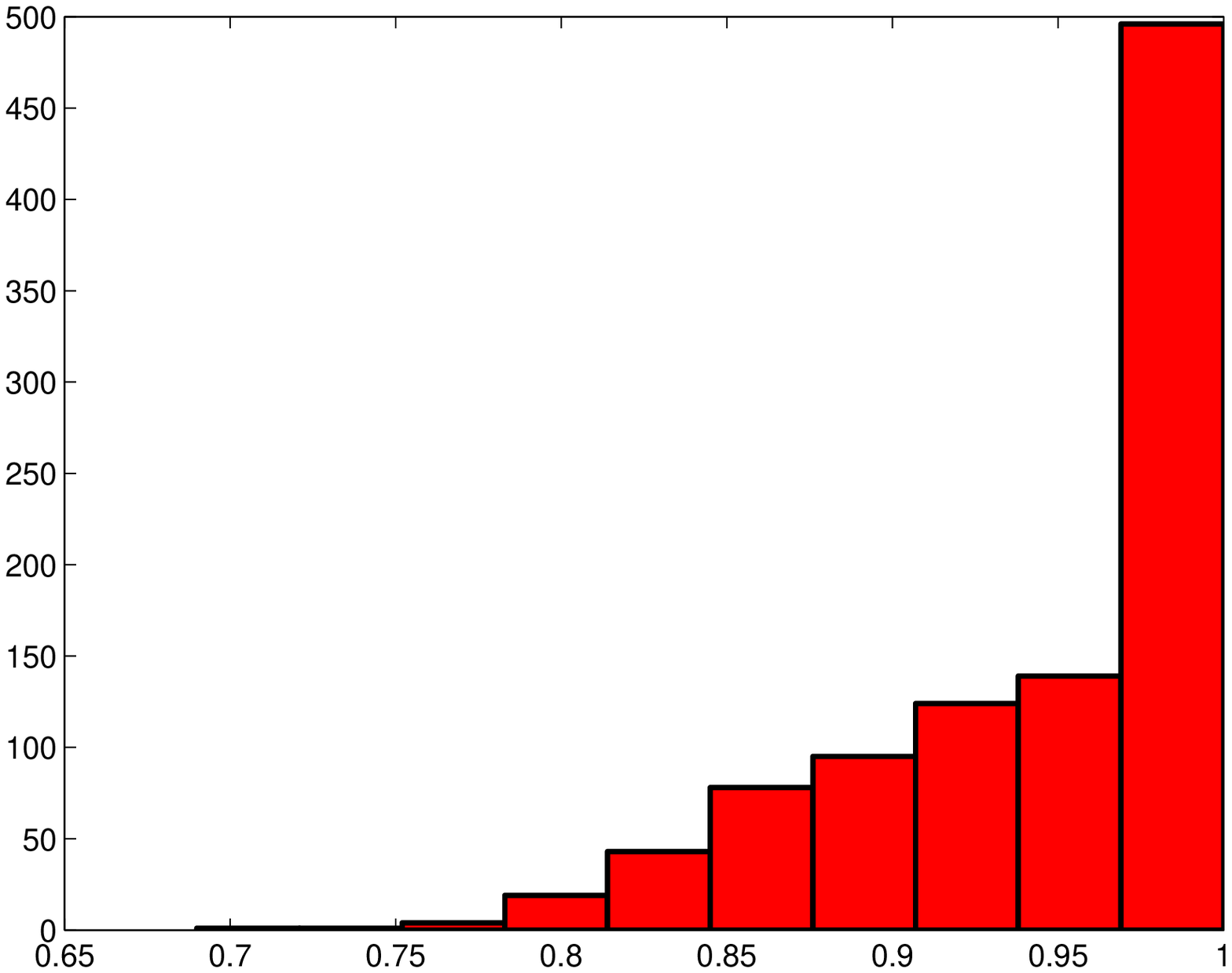} &
    \epsfxsize=2.6in
    \epsffile{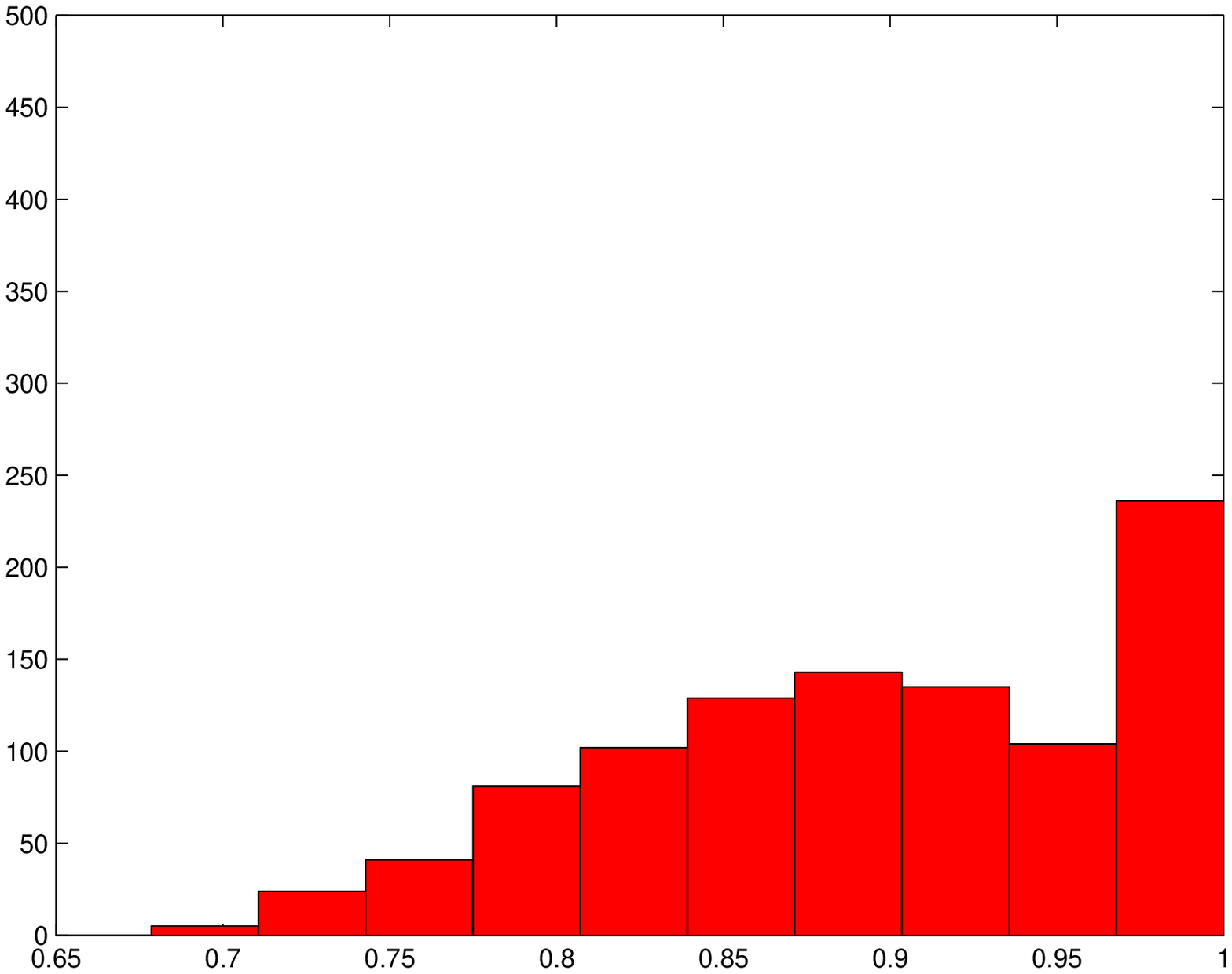}
 \end{array}$
\end{center}
\end{figure}

According to our simulations for $n=10$ in more than 97\% of all instances the copositive bound coincided
with the exact solution, whereas for the Nesterov bound this was true only in 27\% of all cases. There was
not a single instance where the two SDP relaxation bounds would have coincided while being different from
$\alpha ^{\rm ref}$. To assess the quantitative difference between $\alpha _{Q_C} ^{\rm cop}$  and $\alpha
_{C} ^{\rm Nest}$  we show descriptive statistics of the ratio of the two bounds in Table \ref{table1} and we
have plotted  the histogram of the ratio in Figure 1.  Note that for the $3 \times 3$ matrix (\ref{strictex})
with $a=1$ we had a ratio 0.75, which is already quite extreme compared to the minimum 0.69 of our
simulations.

For $20 \times 20$ matrices we obtained similar results. Still in 95\% of all cases the copositive bound was
exact, whereas for the Nesterov bound this ratio dropped to 12\%. There was only one instance for which the
copositive and the Nesterov bound coincided but were smaller than the reference solution. The results of
Table \ref{table1} and the comparison of Figure 1a and Figure 1b indicate that for larger dimensions $n$ the
discrepancy between the two bounds is increasing.

\begin{table}[h]
 \caption{Descriptive statistics for the ratio $ \alpha _{Q_C} ^{\rm cop}/ \alpha _{C} ^{\rm Nest}$  for 1000 randomly generated $n\times n$ matrices.
\label{table1} }
\begin{center}
 \begin{tabular}{|c|c|c|c|c|}
\hline
n& Mean&Std&Min&Median\\
\hline
10& 0.948& 0.058 & 0.690& 0.968\\
\hline
20& 0.896& 0.077 & 0.679 &0.899\\
\hline
\end{tabular}
\end{center}
\end{table}


\section{\label{ellp} Extensions for homogeneous quadratic\\
optimization over the $\ell^p$-ball, $1<p<2$}

We briefly want to address the more general problem
\beq{elpdef} \alpha_p ^* = \min \lk \qf yC  : y\in B_p\rk
\, , \eeq
where 
$1\le p \le 2$. In the case of positive-semidefinite matrices $C$ obviously $\alpha_p
^* = 0$ for all $p$, therefore we concentrate on matrices $C \notin \PP$.

In the SDP Handbook~\cite{NestWolkYe:00}, Nesterov mentions that for $1<p<2$, no practical SDP bounds for~\reff{elpdef}
are in
sight. Because the considered balls $B_p$ are included in the $\ell^2$-ball $B_2$, a cheap bound is always
given by the eigenvalue bound, which is the $\ell^2$-solution
$$\lambda_{\rm min} (C) = \min \lk \qf yC : y\in B_2 \rk\,.$$
Obviously this bound is getting worse for $p$ close to 1. However, we can make use of $\alpha _{Q_C} ^{\rm
cop}$ to obtain a better bound for small $p$\,:

\begin{thm}
A valid SDP-based bound for the problem~\reff{elpdef} is given by
$$\alpha _p
^{\rm cop} := n^{\frac{2(p-1)}{p}}\, \alpha _{Q_C} ^{\rm cop} \ .
$$
\end{thm}
\begin{proof}
We start by blowing up the $\ell^1$-ball such that the result contains $B_p$.  Applying H\"older's inequality
we obtain $\|y\|_1  \le \|y\|_{p}\, \|e\|_{q}= n^{\frac{p-1}{p}} \|y\|_p$ with $\frac 1p + \frac 1q = 1$,
and thus it is evident that $B_p \subseteq n^{\frac{p-1}{p}} B_1$. Therefore
 $$
\alpha_p ^* \ge \min \lk \qf yC : {\norm y }_1 \le n^{\frac{p-1}{p}}\rk \ge n^{\frac{2(p-1)}{p}}\, \alpha
_{Q_C} ^{\rm cop} \, ,
$$
where the last inequality follows from homogeneity of degree two of the quadratic form, and from the validity
of the copositive bound.
\end{proof}

\begin{figure}    
\begin{center}
\caption{Quality of lower bounds in dependence of $p$ for:} $\begin{array}{c@{\hspace{1in}}c}
\multicolumn{1}{l}{{\mbox{\bf (a)}}\quad \mbox{ a typical $10\times 10$ matrix}} &
    \multicolumn{1}{l}{{\mbox{\bf (b)}}  \quad \mbox{the $3\times 3$ matrix}\ C\ \mbox {of \reff{pnormex}}}\\ [-0.03cm]
\epsfxsize=2.6in \epsffile{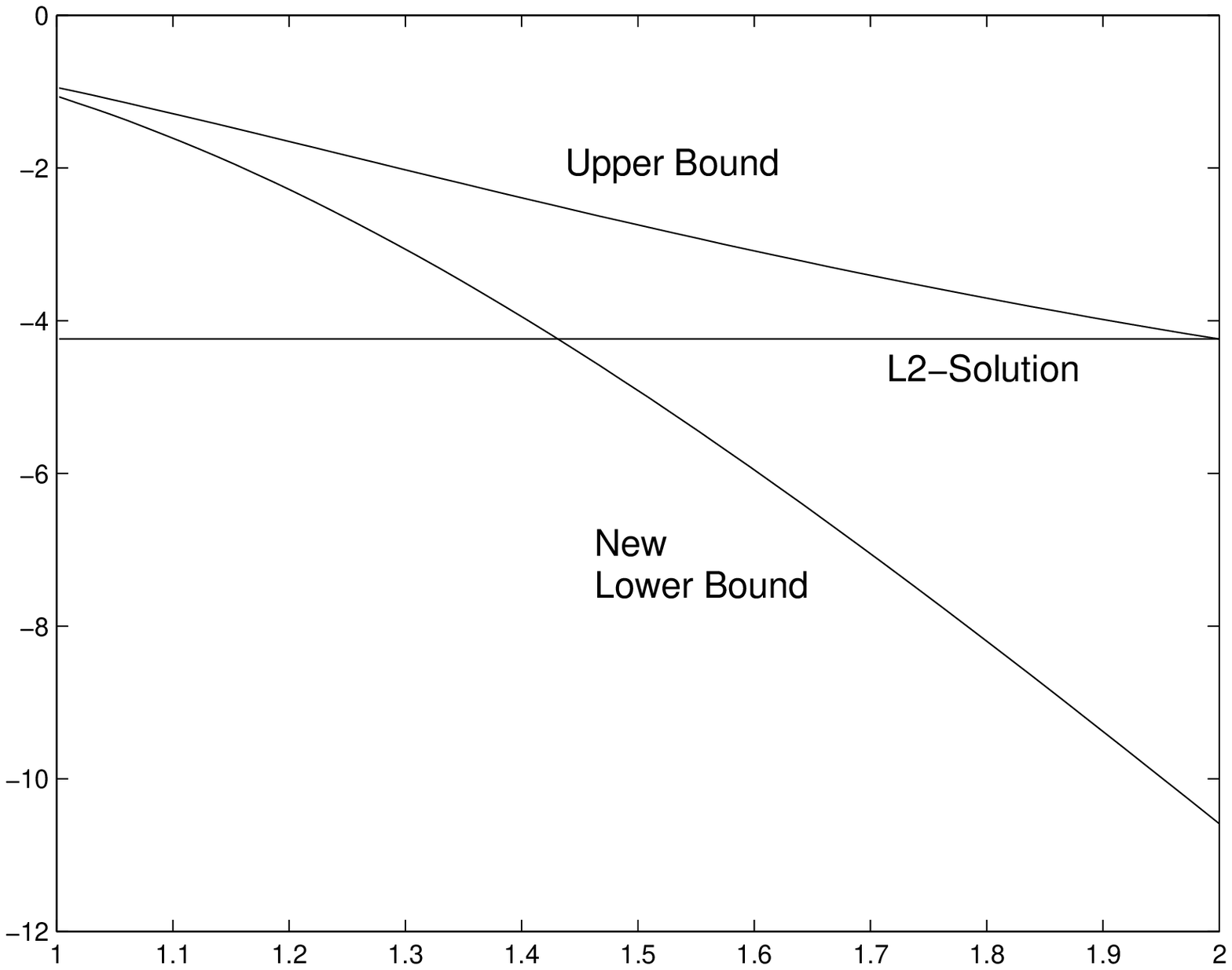} &
    \epsfxsize=2.6in
    \epsffile{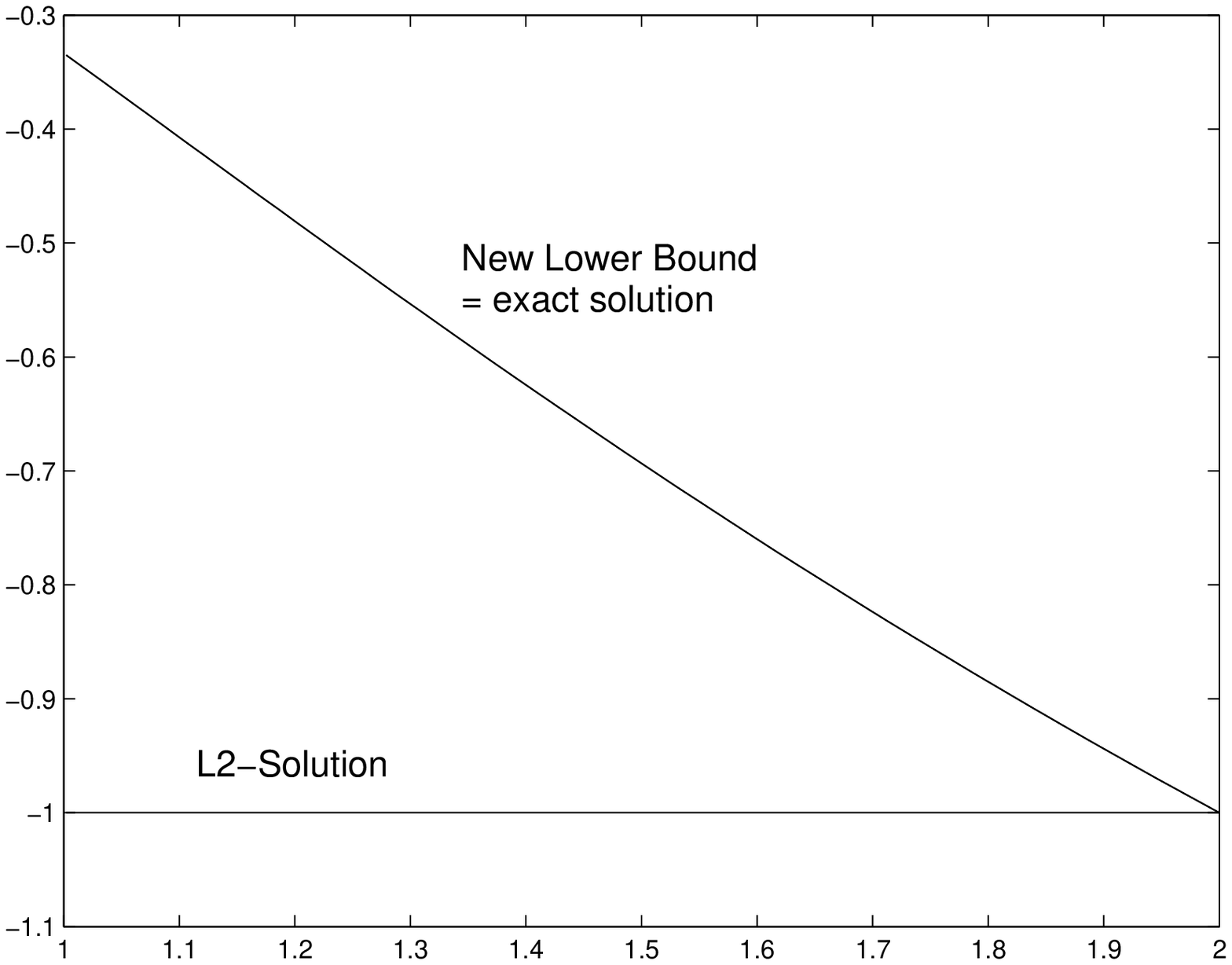}
 \end{array}$
\end{center}
\end{figure}

\noindent We want to discuss the qualities of the lower bound $\alpha _p ^{\rm cop}$ with two examples. In
Figure 2a we consider a randomly generated $10 \times 10$ matrix which illustrates the typical situation: for
$p$ close to one our new bound gives a considerable improvement compared to the eigenvalue bound, for $p=1$
in this case it is actually again the exact solution.  Denoting the ${\norm .}_2$-normalized eigenvector corresponding to
$\lambda_{\rm min} (C)$ by $v_{\rm min}$, we have ${\|v_{\rm min}\|_p}^{-1}v_{\rm min} \in B_p$, hence feasible 
for \reff{elpdef}. Based on that observation we can calculate the simple upper bound
$\frac{\lambda_{\rm min}(C)}{\|v_{\rm min}\|^2_p}$ as plotted in Figure~2a.

In general $\alpha _p ^{\rm cop}$ will be smaller than $\lambda_{\rm min} (C)$ for $p$ close to two. But
there are actually cases where $\alpha _p ^{\rm cop}$ is always larger than $\lambda_{\rm min} (C)$, as in
Figure 2b for the matrix \beq{pnormex}
C=\left[ \begin{array}{rrr} 1 &-1 &1 \\
-1 &1 &1 \\1 &1 &1
\end{array}\right] \, .
\eeq In fact, here $\alpha _p ^{\rm cop}=\alpha_p ^*$ is the exact solution.
The simple upper
bound we considered is given here by $3^{\frac{p-2}{p}} \lambda_{\rm min}(C)$. Now $\alpha _{Q_C} ^{\rm cop}=
-1/3$ and $\lambda_{\rm min}(C)=-1$ leads to equality of the upper and the lower bound.

\end{document}